\documentclass[12pt]{article}

\title{Directional derivatives and subdifferentials of set-valued convex functions}

\author{
Andreas H. Hamel \footnote{Yeshiva University New York, Department of Mathematical
Sciences, \href{mailto:hamel@yu.edu}{hamel@yu.edu}}, Carola Schrage
\footnote{\href{mailto:carolaschrage@gmail.com}{carolaschrage@gmail.com}}
}
\date{{\small \today}}

\usepackage{array} 
\usepackage{verbatim} 

\usepackage{amsmath, amssymb, enumerate}

\newtheorem{theorem}{Theorem}
\newtheorem{corollary}[theorem]{Corollary}
\newtheorem{remark}[theorem]{Remark}
\newtheorem{lemma}[theorem]{Lemma}
\newtheorem{definition}[theorem]{Definition}
\newtheorem{proposition}[theorem]{Proposition}
\newtheorem{example}[theorem]{Example}

\numberwithin{equation}{section}  
\numberwithin{figure}{section}    
\numberwithin{table}{section}     
\numberwithin{theorem}{section}

\newcommand{\of}[1]{\ensuremath{\left( #1 \right)}}

\newcommand{\cb}[1]{\ensuremath{ \left\{ #1 \right\} }}
\newcommand{\sqb}[1]{\ensuremath{ \left[ #1 \right] }}

\newcommand{\bs}{\backslash}

\newcommand{\pend}{ \hfill $\square$ \medskip}

\newcommand{\eps}{\ensuremath{\varepsilon}}

\newcommand{\vp}{\ensuremath{\varphi}}

\renewcommand{\P}{\ensuremath{\mathcal{P}}}

\newcommand{\A}{\ensuremath{\mathcal{A}}}
\newcommand{\F}{\ensuremath{\mathcal{F}}}
\newcommand{\G}{\ensuremath{\mathcal{G}}}

\newcommand{\R}{\mathrm{I\negthinspace R}}
\newcommand{\OLR}{\overline{\mathrm{I\negthinspace R}}}

\newcommand{\dom}{{\rm dom \,}}
\newcommand{\epi}{{\rm epi \,}}

\newcommand{\gr}{{\rm graph \,}}

\newcommand{\cl}{{\rm cl \,}}

\newcommand{\co}{{\rm co \,}}
\newcommand{\core}{{\rm core \,}}

\newcommand{\isum}{{+^{\negmedspace\centerdot\,}}}

\newcommand{\idif}{{-^{\negmedspace\centerdot\,}}}

\newcommand{\triup}{{\rm \vartriangle}}
\newcommand{\trido}{{\rm \triangledown}}

\newcommand{\cone}{{\rm cone\,}}
\newcommand{\Int}{{\rm int\,}}

\newcommand{\istardif}{-_{z^*}}

\usepackage[
colorlinks=true, linkcolor=blue, citecolor=blue, urlcolor=blue, filecolor=blue
auskommentieren und unten pdfborder aktivieren]{hyperref}

\begin{document}
\maketitle

\begin{abstract}
A new directional derivative and a new subdifferential for set-valued convex functions are constructed, and a set-valued version of the so-called 'max-formula' is proven. The new concepts are used to characterize solutions of convex optimization problems with a set-valued objective. As a major tool, a residuation operation is used which acts in a space of closed convex, but not necessarily bounded subsets of a topological linear space. The residuation serves as a substitute for the inverse addition and is intimately related to the Minkowski or geometric difference of convex sets. The results, when specialized, even extend those for extended real-valued convex functions since the improper case is included.
\end{abstract}

\pagebreak

\section{Introduction}

In this note, we introduce new notions of directional derivatives and subdifferentials for set-valued convex functions, we prove the so-called max-formula, a result of "exceptional importance" in the scalar case \cite[p. 90]{Zalinescu02}, and characterize solutions of set-valued optimization problems in terms of the new derivatives. The latter topic sheds some new light on what should actually be understood by a solution of a convex optimization problem with a set-valued objective. In particular, we supplement the solution concept given in \cite{HeydeLoehne11} by a new one and show that a solution set can be reduced to a singleton via a generalized translation.

There exists basically three different (but partially overlapping) approaches for defining derivatives for set-valued functions. One approach starts by picking a point in the graph of the set-valued function and assigns to it another set-valued function whose graph is  some kind of tangent cone to the graph of the original function at the point in question. The book \cite{AubinFrankowska90Book} gives a prestigious account of such concepts, and Mordukhovich's coderivative \cite{Mordukhovich06-1Book} is of the same nature. The second approach selects a class of 'simple' set-valued functions the elements of which shall serve as approximation for a general set-valued function and then defines what is actually understood by "approximation". A representative for this approach is \cite{LemarechalZowe91}.
The third approach embeds the class of set-valued functions under consideration into a linear space and operates with classical derivative concepts. The reader may consult \cite{GorokhovikZabreiko05} for more references and a more complete account of the three basic approaches described above. Note, however, that the last two approaches very often are restricted to set-valued functions with compact convex values, and to finite dimensional \cite{GorokhovikZabreiko05} or even one-dimensional pre-image spaces \cite{LemarechalZowe91}.

On the other hand, it turned out that it is a hard task to generalize basic results in convex analysis from extended real- to vector- or even set-valued functions. The 'max-formula' may feature as an example which is relevant for the present paper: Under some qualifying conditions, the directional derivative of a convex function at a given point is the support function of the subdifferential at the same point. Since this implies the non-emptiness of the subdifferential, this result is counted among the 'core results of the convex analysis' \cite[p. 122]{BorweinZhu05}. The difficulties which arise when passing from one-dimensional to more general  image spaces are brought out, for example, in \cite[Theorem 6.1]{Borwein82}: The pre-image space must be a "Minkowski differentiability space", the image space must be ordered by a closed normal cone and enjoy the so-called monotone sequence (= greatest lower bound) property.

Our approach is more in the spirit of traditional derivative concepts which rely on increments of a function at a point in some direction. Using a residuation instead of a difference (an inverse group operation which is not available in relevant subsets of the power set of a linear space) we are able to define difference quotients and their limits even for set-valued functions. In fact, it seems to be natural to "skip" the vector-valued case by embeding it into the set-valued one. The residuation is defined on carefully selected subsets of the power set of the (linear) image space; these subsets carry the order structure of a complete lattice (= every subset has an infimum and supremum) and the algebraic structure of a semi-module over the semi-ring $\R_+$. It turns out that the old concept of the Minkowski (or geometric) difference of convex sets \cite{Hadwiger50} can be identified with the residuation in these spaces of sets; even this seems to be a new contribution (see also \cite{HamelSchrage12}) although residuations have been used before in (convex) analysis, see for example \cite{GetanMartinezLegazSinger03}, \cite{CohenGaubertQuadratSinger05} and the references therein.

The dual variables in our theory are simple set-valued functions generated by pairs of continuous linear functionals instead of continuous linear operators as, for instance, in \cite{Borwein82}, \cite{Borwein84}. Moreover, no restrictive assumptions to the ordering cone in the underlying (linear) image space are imposed such as normality, pointedness, non-empty interior, generating a lattice order etc. These features make the theory presented in this note much more adequate for applications. The interested reader is referred to \cite{HamelHeydeRudloff11} for a financial application where the ordering cone is not pointed, in general, and has 'many' generating vectors.

In the next section, the basics about set-valued functions and their image spaces are introduced. Section \ref{SecDirDer} contains the definitions of directional derivatives and subdifferentials for set-valued convex functions and the main results. Section \ref{SecConProc} presents the link between 'adjoint process duality' (Borwein, 1983), Mordukhovich's coderivative and our derivative concepts. In the final section, set-valued optimization problems are discussed.

\section{Preliminaries}

\subsection{Image spaces}

Let $Z$ be a locally convex, topological linear space and $C \subseteq Z$ a convex cone with $0 \in C$.
We write $z_1\leq_C z_2$ for $z_2 - z_1 \in C$ with $z_1, z_2 \in Z$ which defines a
reflexive and transitive relation (a preorder). The topological dual space
of $Z$ is denoted by $Z^*$, the (positive) dual cone of $C$ by $C^+ = \cb{z^* \in Z^*
\;\vert\; \forall z \in C \colon  z^*\of{z} \geq 0}$. Note that $C^+ \neq \cb{0}$ if, and only if, $\cl C \neq Z$ which is assumed throughout the paper. The negative dual cone is $C^-=-C^+$.

The relation $\leq_C$ on $Z$ can be extended to the powerset $\P\of{Z}$ of $Z$, the set of
all subsets of $Z$ including the empty set $\emptyset$ in two canonical ways (see \cite{Hamel09}
and the references therein). This gives rise to consider the following subsets of
$\P\of{Z}$:
\begin{align*}
\F\of{Z, C} & = \cb{A \in \P\of{Z} \;\vert\; A = \cl\of{A+C}} \\
 \label{EqDefQ} \G\of{Z, C} & = \cb{A \in \P\of{Z} \;\vert\; A = \cl\co\of{A+C}}.
\end{align*}
Elements of $\F\of{Z, C}$ are sometimes called upper closed (\cite[Definition 1.50]{Loehne11Book}) with respect to $C$. We shall abbreviate $\F\of{Z, C}$ and $\G\of{Z, C}$ to $\F\of{C}$ and $\G\of{C}$, respectively.

The Minkowski (elementwise) addition for non-empty subsets of $Z$ is extended to
$\P\of{Z}$ by
\[
\emptyset + A = A + \emptyset = \emptyset
\]
for $A \in \P\of{Z}$. Using this, we define an associative and commutative binary operation $\oplus \colon
\F\of{C} \times \F\of{C} \to \F\of{C}$ by
\begin{eqnarray}
\label{EqSumInP} A \oplus  B= \cl\of{A+B}
\end{eqnarray}
for $A, B \in \F\of{C}$. The elementwise multiplication of a set $A \subseteq Z$ with a
(non-negative) real number is extended by
\[
0 \cdot A = \cl C, \quad t \cdot \emptyset = \emptyset
\]
for all $A \in \F\of{C}$ and $t > 0$. In particular, $0 \cdot \emptyset = \cl C$ by
definition, and we will drop the $\cdot$ in most cases.

The triple $\of{\F\of{C}, \oplus, \cdot}$ is a conlinear space with neutral element $\cl C$,
and, obviously, $\of{\G\of{C}, \oplus, \cdot}$ is a conlinear subspace of it. The concept of a 'conlinear space' has been introduced in \cite{Hamel05Habil}, see also \cite{Hamel09}, \cite{HamelSchrage12}. It basically means that $\of{\F\of{C}, \oplus}$ is a commutative monoid, and a multiplication of elements of $\F\of{C}$ with those of $\R_+$ is defined and satisfies some obvious requirements, but not, in general, the second distributivity law $\of{s + t} \cdot A = s \cdot A \oplus t \cdot A$ for $s, t \in \R_+$, $A \in \F\of{C}$. The elements of $\F\of{C}$ which do satisfy this law are precisely those of $\G\of{C}$, thus $\of{\G\of{C}, \oplus, \cdot}$ is a semi-module over the semi-ring $\R_+$.

On $\F\of{C}$ and $\G\of{C}$, $\supseteq$ is a partial order which is compatible with the
algebraic operations just introduced. Thus, $\of{\F\of{C}, \oplus, \cdot, \supseteq}$ and
$\of{\G\of{C}, \oplus, \cdot, \supseteq}$ are partially ordered, conlinear spaces in the sense
of \cite{Hamel05Habil}, \cite{Hamel09}. Note that this is true without any further assumptions to $C$.
In particular, $C$ is not required to generate a partial order, a fact, which will be
used later on.

We will abbreviate $\F^\triup = \of{\F\of{C}, \oplus, \cdot, \supseteq}$ and
$\G^\triup = \of{\G\of{C}, \oplus, \cdot, \supseteq}$, and we will write $A \in \G^\triup$ and $\mathcal A \subseteq \G^\triup$ in order to denote an element $A \in \G\of{C}$ and a subset $\mathcal A \subseteq \G\of{C}$, respectively.

Moreover, $\of{\F\of{C}, \supseteq}$ and $\of{\G\of{C}, \supseteq}$
are complete lattices with greatest (top) element $\emptyset$ and least (bottom) element
$Z$. For a subset $\A \subseteq \G^\triup$, the infimum and the supremum of $\A$ are
given by
\begin{eqnarray}
\label{EqInf} \inf \A = \cl\co \bigcup\limits_{A \in \A} A,\qquad
\sup \A = \bigcap\limits_{A \in \A} A
\end{eqnarray}
where we agree upon $\inf \A = \emptyset$ and $\sup \A = Z$
whenever $\A = \emptyset$. Finally, for all $\A \subseteq \G^\triup$ and $B \in
\G^\triup$,
\begin{eqnarray}
\label{EqResiduEq} B \oplus \inf\A = \inf\of{B \oplus \A}
\end{eqnarray}
where $B \oplus \A = \cb{B \oplus A \mid A \in \A}$. It follows that $\G^\triup$ is an
$\inf$-residuated space (see \cite{HamelSchrage12} for more details). The inf-residuation will
serve as a substitute for the inverse addition and is defined as follows: For $A, B\in
\G^\triup$, set
\begin{equation}
\label{EqInfResiduation} A \idif B = \inf\cb{D \in \G^\triup \mid B + D \subseteq A} =
\cb{z\in Z \mid B + z\subseteq A}.
\end{equation}
Note that, for $A \in \mathcal G^\triup$, the set on the right hand side of \eqref{EqInfResiduation} is indeed closed since
\[
\cb{z \in Z \;\vert\; B + z\subseteq A} = \bigcap_{b \in B}\cb{z \in Z \;\vert\; b + z \in A}
\]
which is an intersection of closed sets whenever $A$ is closed.

Sometimes, the right hand side of \eqref{EqInfResiduation} is called the geometric difference \cite{Pontrjagin80} or the Minkowski difference \cite{Hadwiger50} of the two sets $A$ and $B$, and H. Hadwiger should probably be credited for its introduction. The relationship with residuation theory  (see, for instance, \cite{BlythJanowitz72Book}, \cite{GalatosJipsenKowalskiOno07}) has been established in \cite{HamelSchrage12}. At least, we do not know an earlier reference. 

\begin{example}
\label{ExExtReals}
Let us consider $Z = \R$, $C = \R_+$. Then $\G\of{Z, C} = \cb{[r, +\infty) \mid r \in \R}\cup\cb{\R}\cup\cb{\emptyset}$, and $\G^\triup$ can be identified (with respect to the algebraic and order structures which turn $\G\of{\R, \R_+}$ into an ordered conlinear space and a  complete lattice admitting an inf-residuation) with $\OLR = \R\cup\cb{\pm\infty}$ using the 'inf-addition' $\isum$ (see \cite{RockafellarWets98} , \cite{HamelSchrage12}). The inf-residuation on $\OLR$ is given by
\[
r \idif s  = \inf\cb{t\in\R \mid r \leq s \isum t}
\]
for all $r,s\in\OLR$, compare \cite{HamelSchrage12} for further details.
\end{example}

Historically, it is interesting to note that R. Dedekind \cite{Dedekind1872} introduced the residuation concept and used it in order to construct the real numbers as 'Dedekind sections' of rational numbers. The construction above is in this line of ideas, but in a rather abstract setting.

\begin{remark}
The inf-residuation can be defined on $\mathcal F^\triup$ and even other subspaces of $\mathcal P\of{Z}$, but we only need the construction in $\mathcal G^\triup$ in this paper. Likewise, in $\G^\trido = \of{\G\of{-C}, \oplus, \cdot, \subseteq}$, a sup-residuation can be defined such that  the whole theory becomes symmetric. The interested reader is referred to \cite{HamelSchrage12}. 
\end{remark}

In many cases, the set $A \idif B$ is "too small", even empty: Consider $Z = \R^2$, $A = C = \R^2_+$, $B = \cb{z \in \R^2 \mid \eps z_1 + z_2 \geq 0, \; z_1 + \eps z_2 \geq 0}$. Then $A \idif B = \emptyset$ for each $\eps > 0$. Therefore, we modify the inf-residuation in $\G^\triup$ as follows. Take $z^* \in C^-\bs\cb{0}$ and let $H\of{z^*} = \cb{z \in Z \mid z^*\of{z} \leq 0}$ be the homogenous half-space with normal $z^*$. We set
\begin{eqnarray}
\label{EqZstarResi} A \istardif B = \of{A \oplus H\of{z^*}} \idif B =
    \cb{z \in Z \mid B + z\subseteq A \oplus H\of{z^*}}.
\end{eqnarray}

The operation $\istardif$ can be expressed using the inf-residuation in $\OLR$ and support functions, see \cite[Proposition 5.20]{HamelSchrage12} and therefore, it would be interesting to study the relationships to the Demyanov and Rubinov difference \cite[p. 180 and p. 182, respectively]{RubinovAkhunov92Opt}. However, our construction is particularly tailored for non-compact convex sets.

By definition, $A \istardif B = Z$ if $A \oplus H\of{z^*} = Z$ or $B = \emptyset$, and $A \istardif B = \emptyset$ if $A \oplus H\of{z^*} \neq Z$, $B \oplus H\of{z^*} = Z$ and if $A = \emptyset$, $B \neq \emptyset$. In all other cases, $A \istardif B$ is a non-empty closed half-space parallel to $H\of{z^*}$. The relationships
\[
A \istardif B = \of{A \oplus H\of{z^*}} \istardif B = A \istardif \of{B \oplus H\of{z^*}} = 
	\of{A \oplus H\of{z^*}} \istardif \of{B \oplus H\of{z^*}}.
\] 
and 

\begin{equation}\label{EqAMinusPlusB}
 A \oplus H\of{z^*}  \supseteq B \oplus \of{A\istardif B}
\end{equation}
for all $A, B \in \G^\triup$ are immediate from the definition of $\istardif$. The next proposition makes it clear that the expression $H\of{z^*} \istardif A$ replaces
$-A$.

\begin{proposition}
\label{MinusSetInequality} Let $A, B \in \G^\triup$ and $z^* \in C^-\bs\cb{0}$. Then (a)
\[
A \oplus H\of{z^*} \supseteq B \oplus H\of{z^*} \quad \Leftrightarrow \quad
    H\of{z^*} \istardif B \supseteq H\of{z^*} \istardif A,
\]
and (b)
\[
A \istardif B \not\in \cb{Z, \emptyset} \quad \Leftrightarrow \quad A \oplus H\of{z^*}, B \oplus H\of{z^*} \not\in \cb{Z, \emptyset}.
\]
\end{proposition}

{\sc Proof.} (a) "$\Rightarrow$": We have
\begin{multline*}
H\of{z^*} \istardif A = H\of{z^*} \istardif \of{A \oplus H\of{z^*}} 
	= \cb{z \in Z \mid  A \oplus H\of{z^*} + z \subseteq H\of{z^*}} \\
\subseteq \cb{z \in Z \mid  B \oplus H\of{z^*} + z \subseteq H\of{z^*}} = H\of{z^*} \istardif B
\end{multline*}
since $B \oplus H\of{z^*} \subseteq A \oplus H\of{z^*}$.

"$\Leftarrow$": This implication is certainly true if $A \oplus H\of{z^*} = Z$. If $A
\oplus H\of{z^*} = \emptyset$, then $H\of{z^*} \istardif  A = Z$, hence $H\of{z^*}
\istardif B = Z$ by assumption which in turn implies $B = \emptyset$. Finally, assume $A
\oplus H\of{z^*} = z_A +  H\of{z^*}$ for some $z_A \in Z$. Then
\[
-z_A \in H\of{z^*} \istardif A \supseteq H\of{z^*} \istardif B,
\]
hence $B \oplus H\of{z^*} \subseteq z_A +  H\of{z^*} \subseteq A \oplus  H\of{z^*}$.

(b) is a straightforward consequence of the definition of $\istardif$.  \pend


\medskip The following calculus rules for $\istardif$ apply and will be used frequently.

\begin{proposition}
\label{PropCalcIstardif} Let $A, B, D \subseteq \G^\triup$ and $z^* \in C^-\bs\cb{0}$. Then

(a) $A \supseteq B$ $\Rightarrow$ $A \istardif D \supseteq
B \istardif D$ and $D \istardif B \supseteq D \istardif A$.

(b) $A \istardif A = H\of{z^*}$ if $A \oplus H\of{z^*} \not\in \cb{Z, \emptyset}$, and $A \istardif A = Z$
if $A \oplus H\of{z^*} \in \cb{Z, \emptyset}$.

(c) $s_1, s_2 \geq 0$, $s_1 + s_2 = 1$ $\Rightarrow$
\[
\of{s_1 A \oplus s_2 B} \istardif D \supseteq s_1\of{A \istardif D} \oplus
s_2\of{B \istardif D}.
\]

(d) $\of{A \oplus B} \istardif B \supseteq A \oplus H\of{z^*}$. The strict inclusion applies if, and only if, $A = B = \emptyset$, or $A \oplus H\of{z^*} \not\in \cb{Z, \emptyset}$ and $B \oplus H\of{z^*} \in \cb{Z, \emptyset}$.

(e) $A \oplus H\of{z^*} \supseteq \of{A \istardif B} \oplus B$. The strict inclusion applies if, and only if, $A \neq B = \emptyset$, or $A \oplus H\of{z^*} \not\in \cb{Z, \emptyset}$ and $B \oplus H\of{z^*} = Z$.

(f) $\of{A \istardif D} \oplus B \subseteq \of{A \oplus B} \istardif D$.  The strict inclusion applies if, and only if, $B = D = \emptyset$, or $B \oplus H\of{z^*} = D \oplus H\of{z^*}  = Z$ and $A \oplus H\of{z^*} \not\in \cb{Z, \emptyset}$. 

(g) $H\of{z^*} \supseteq A \istardif B$ $\Rightarrow$ $A \subseteq B \oplus H\of{z^*}$.
\end{proposition}

{\sc Proof.} (a) - (c) are elementary using the definition of $\istardif$.

(d) The inclusion immediately follows from the definition of $\istardif$. If $A \oplus H\of{z^*} \not\in \cb{Z, \emptyset}$, $B
\oplus H\of{z^*} \not\in \cb{Z, \emptyset}$ we can find $z_A, z_B \in Z$ such
that
\[
A \oplus H\of{z^*} = z_A + H\of{z^*}, \quad B \oplus H\of{z^*} = z_B + H\of{z^*}.
\]
Then
\[
B \oplus H\of{z^*} + \of{z_A - z_B} = z_B + \of{z_A - z_B} + H\of{z^*}
    = z_A + H\of{z^*} = A \oplus H\of{z^*},
\]
hence $z_A - z_B \in A \istardif B$ and $z_A - z_B + H\of{z^*} = A \istardif B$. This
implies
\[
\of{z_A - z_B + H\of{z^*}} \oplus B = \of{A \istardif B} \oplus B = A \oplus H\of{z^*}.
\]
In view of proposition \ref{MinusSetInequality} (b), this leaves two possibilities for the strict inclusion: The first is $A = \emptyset$ and $\of{A \oplus B} \istardif B \neq \emptyset$, the second $A \oplus H\of{z^*} \not\in \cb{Z, \emptyset}$ and $\of{A \oplus B} \istardif B = Z$. In the first, we obtain $\of{A \oplus B} \istardif B = \of{\emptyset \oplus B} \istardif B = \emptyset \istardif B$. The set $\emptyset \istardif B$ is non-empty if, and only if, $B = \emptyset$ in which case $\emptyset \istardif B = Z$, so strict inclusion holds. In the second case, $\of{A \oplus B} \istardif B = Z$ precisely when $B \oplus H\of{z^*} \in \cb{Z, \emptyset}$. Together, we obtain the conditions in (d).

(e) This claim can be proven by similar arguments as used for (d).

(f)  If $\of{A \istardif D} \oplus B = \emptyset$ then the inclusion is trivially true. Otherwise, for each $z \in \of{A \istardif D}$ (there is one!) $D + z \subseteq A \oplus H\of{z^*}$ which implies $D + B + z \subseteq A \oplus B \oplus H\of{z^*}$ which in turn
gives $B + z \subseteq \of{A \oplus B} \istardif D$. The inclusion follows. 

If both sides are neither $Z$ nor $\emptyset$, then equality holds true. Indeed, in this case (see proposition \ref{PropCalcIstardif} (b)) $A \oplus H\of{z^*}, B \oplus H\of{z^*}, D \oplus H\of{z^*} \not\in \cb{Z, \emptyset}$, hence there are $z_A, z_B \in Z$ such that $A \oplus H\of{z^*} = z_A + H\of{z^*}$ and $B \oplus H\of{z^*} = z_B + H\of{z^*}$. This gives
\begin{align*}
\of{A \istardif D} \oplus B
    & = \cl\cb{z + b \;\vert\; z \in Z, \, b \in B, \,
        D + z \subseteq z_A + H\of{z^*}} \\
    & = \cb{z + z_B \;\vert\; z \in Z, \,
        D + z \subseteq z_A + H\of{z^*}} \oplus H\of{z^*} \\
    & = \cb{z \;\vert\; z \in Z, \,
        D + z - z_B \subseteq z_A + H\of{z^*}} = \of{A \oplus B} \istardif D.
\end{align*}
This leaves two cases for strict inclusion: The first is $\of{A \istardif D} \oplus B = \emptyset$, the second $\of{A \istardif D} \oplus B \not\in\cb{Z, \emptyset}$ and $\of{A \oplus B} \istardif D = Z$. The second case can not occur as a little straightforward analysis shows. In the first, we can have $B = \emptyset$ in which case $\of{A \oplus B} \istardif D = \emptyset$ and $ \of{A \oplus B} \istardif D \neq \emptyset$ if, and only if, $D = \emptyset$. Or we can have $A \istardif D = \emptyset$ which produces the strict inclusion precisely when  $B \oplus H\of{z^*} = D \oplus H\of{z^*}  = Z$ and $A \oplus H\of{z^*} \not\in \cb{Z, \emptyset}$.

(g) First, note that the assumption implies $A \oplus H\of{z^*} \neq Z$. If $B \oplus
H\of{z^*} = Z$, $A \oplus H\of{z^*} \neq Z$ then $A \istardif B = \emptyset$, the
assumption is trivially satisfied, and, also trivially, $A \subseteq B \oplus H\of{z^*}$.

Now, assume $A \oplus H\of{z^*} \neq Z$, $B \oplus H\of{z^*} \neq Z$ and that $A
\subseteq B \oplus H\of{z^*}$ is not true. Then $B \subseteq \Int\of{A \oplus H\of{z^*}}
= A + \Int H\of{z^*}$. Therefore, we can find $z \in Z$ such that $z^*\of{z} >0$ and $B+z
\subseteq A \oplus H\of{z^*}$ contradicting the assumption. \pend

\begin{lemma} Let $A, B \in \G^\triup$ and $z^* \in C^-\bs\cb{0}$. Then
\begin{equation}\label{EqAMinusBminusA}
H\of{z^*} \istardif B \subseteq \of{A \istardif B} \istardif A.
\end{equation}
The strict inclusion applies if, and only if, $B \oplus H\of{z^*} = Z$ and $A \oplus H\of{z^*} \in \cb{Z, \emptyset}$, or $B \oplus H\of{z^*}  \neq \emptyset$ and $A \oplus H\of{z^*} = Z$.
\end{lemma}

{\sc Proof.}  By proposition \ref{PropCalcIstardif} (e) we have
\[
B \oplus  \sqb{H\of{z^*} \istardif B}  \subseteq H\of{z^*}, 
\]
hence
\[
A \oplus B \oplus  \sqb{H\of{z^*} \istardif B}  \subseteq A \oplus H\of{z^*}.
\]
The latter is equivalent to \eqref{EqAMinusBminusA}. Indeed, by definition of $\istardif$, \eqref{EqAMinusBminusA} is equivalent to 
\[
A \oplus \sqb{H\of{z^*} \istardif B}  \subseteq A \istardif B = \sqb{A \oplus H\of{z^*}} \istardif B
\]
which in turn is equivalent to 
\[
B \oplus \sqb{A \oplus \sqb{H\of{z^*} \istardif B}} \subseteq A \oplus H\of{z^*}.
\]
Altogether, this proves \eqref{EqAMinusBminusA}. 

By proposition \ref{MinusSetInequality} (b), $\of{A \istardif B} \istardif A \not\in \cb{Z, \emptyset}$ if, and only if, $A \oplus H\of{z^*} = z_A + H\of{z^*}$, $B \oplus H\of{z^*} = z_B + H\of{z^*}$ for some $z_A, z_B \in Z$. In this case, \eqref{EqAMinusBminusA} is satisfied as an equation. Therefore, the strict inclusion applies if, and only if, $H\of{z^*} \istardif B \neq Z$ and $\of{A \istardif B} \istardif A = Z$. One directly checks that this amounts to $B \oplus H\of{z^*} = Z$ and $A \oplus H\of{z^*} \in \cb{Z, \emptyset}$, or $B \oplus H\of{z^*}  \neq\emptyset$ and $A \oplus H\of{z^*} = Z$.  \pend

\subsection{$\G^\triup$-valued functions and their scalarizations}

Let $X$ be another locally convex, topological linear space with dual $X^*$. A function $f \colon X \to
\G^\triup$ is called convex if
\begin{equation}
\label{EqConvFct} \forall x_1, x_2 \in X, \; \forall t \in \of{0,1} \colon
 f\of{tx_1+(1-t)x_2} \supseteq tf\of{x_1} \oplus \of{1-t}f\of{x_2}.
\end{equation}
It is an exercise (see, for instance, \cite{Hamel09}) to show that $f$ is convex if, and only if, the set
\[
\gr f = \cb{\of{x,z} \in X \times Z \colon z \in f\of{x}}
\]
is convex. A $\G^\triup$-valued function $f$ is called positively homogeneous if
\[
\forall t > 0, \forall x \in X \colon f\of{tx} \supseteq tf\of{x},
\]
and it is called sublinear if it is positively homogeneous and convex. Another exercise shows that $f$ is sublinear if, and only if, $\gr f$ is a convex cone.

A function $f \colon X \to \G^\triup$ is called lower semi-continuous (l.s.c. for short)
at $x_0 \in X$ if $f\of{x_0} \supseteq \liminf\limits_{x \to x_0}f\of{x}$ where
\begin{eqnarray}
\label{EqLimInf} \liminf\limits_{x \to x_0}f\of{x_0} =
 \bigcap\limits_{U \in \mathcal N_X\of{0}}\cl\co\bigcup\limits_{x \in x_0 + U}f\of{x}
\end{eqnarray}
where $\mathcal N_X\of{0}$ is a neighborhood base of $0 \in X$. The function $f$ is
called closed if it is l.s.c. at every $x \in X$. Again, one can show that $f$ is closed
if, and only if, $\gr f \subseteq X \times Z$ is a closed set with respect to the product
topology, see \cite[Proposition 2.34]{Loehne11Book}.

The greatest closed convex minorant of a function $f \colon X \to \G^\triup$ is denoted
by $\of{\cl\co f} \colon X \to \G^\triup$. We have
\begin{eqnarray}
\label{EqClCoHullFunction} \gr\of{\cl\co f}=\cl\co\of{\gr f}.
\end{eqnarray}

\begin{remark} A more common convexity concept for functions $F \colon X \to \P\of{Z}$ is the following (compare, for instance, \cite[Definition 14.6]{Jahn04}): $F$ is called $C$-convex if
\[
t \in \sqb{0,1}, x_1, x_2 \in X \; \Rightarrow \; F\of{tx_1 + \of{1-t}x_2} + C \supseteq tF\of{x_1} +  \of{1-t}F\of{x_2}.
\]
It is easily seen that $C$-convexity of $F$ implies that the function $x \mapsto f\of{x} = F\of{x} + C$ maps into 
\[
\cb{D \in \P\of{Z} \mid D = \co\of{D+C}}
\]
and has a convex graph which coincides with 
\[
\epi F = \cb{\of{x, z} \in X \times Z \mid z \in F\of{x} + C}
\]
(see \cite[Definition 14.7]{Jahn04}). Moreover, if $\gr f$ is additionally closed, then $f$ automatically maps into $\G\of{C}$. 

Finally, note that it does not make sense to distinguish between the graph and the
epigraph of a $\G^\triup$-valued function since the two sets coincide.
\end{remark}

A function $f \colon X \to \G^\triup$ is called proper, if its domain
\[
\dom f = \cb{x \in X \colon f\of{x} \neq \emptyset}
\]
is nonempty and $f$ does not attain the value $Z$. A $\G^\triup$-valued function is called $C$-proper if $\of{f\of{x} - C}\bs f\of{x} \neq \emptyset$ for all $x \in \dom f$. A function is called $z^*$-proper for $z^* \in C^-\bs\cb{0}$ if the function $x \mapsto f\of{x} \oplus H\of{z^*}$ is proper. Of course, if $f$ is $z^*$-proper for at least one $z^* \in C^-\bs\cb{0}$, then it is proper. Vice versa, if $f$ is closed convex proper function, then there is at least one $z^* \in C^-\bs\cb{0}$ such that $f$ is $z^*$-proper. The latter fact follows, for example, from \cite[Theorem 1]{Hamel09}.

\begin{example}
\label{ExConlinFunction}
Let $x^* \in X^*$ and $z^* \in Z^*$ be given. The function $S_{\of{x^*, z^*}} \colon X \to \mathcal P\of{Z}$ defined through
\[
S_{\of{x^*, z^*}}\of{x} = \cb{z \in Z \mid x^*\of{x} + z^*\of{z} \leq 0}
\] 
maps into $\G\of{C}$ if, and only if, $z^* \in C^-$. Moreover, it is positively homogeneous and additive. Therefore, if $z^* \in C^-$ the function  $x \mapsto S_{\of{x^*, z^*}}\of{x}$ is $\G\of{C}$-valued and  convex. It is $sz^*$-proper for all $s>0$ if, and only if, $z^* \in C^-\bs\cb{0}$. Finally, $S_{\of{x^*, z^*}}\of{0} = H\of{z^*} = \cb{z \in Z \mid z^*\of{z} \leq 0}$ is a homogeneous closed half space with normal $z^*$ if $z^* \neq 0$ and $S_{\of{x^*, 0}}\of{0} = Z$. In particular, $S_{\of{x^*, z^*}}\of{0} \neq \cb{0}$ for all $z^* \in C^-$. 

For $z^*\in C^-$, the useful relation
\[
\forall x \in X, \; \forall t> 0 \colon S_{\of{x^*,tz^*}}\of{x} = S_{\of{\frac{1}{t}x^*,z^*}}\of{x}
\]
immediately follows from the definition of $S_{\of{x^*,z^*}}\of{x}$. If $Z=\R$ and $C = R_+$, then
$S_{\of{x^*, -s}}\of{x} = -\frac{1}{s}x^*\of{x} + \R_+$ for $s \in C^-\bs\cb{0} = -\R_+\bs\cb{0}$ while $S_{\of{x^*,0}}\of{x} \in \cb{Z, \emptyset}$.
\end{example}

Let a function $f \colon X \to \G^\triup$ be given. The family of extended real-valued functions $\vp_{f, z^*} \colon  X \to \R\cup\cb{\pm\infty}$ defined by
\[
\vp_{f, z^*}\of{x} = \inf\cb{-z^*\of{z} \mid z \in f\of{x}}, \; z^* \in C^-\bs\cb{0},
\]
is called the family of (linear) scalarizations for $f$. The function $f$ is convex if, and only if, the scalarizing function $\vp_{f, z^*}$ is convex for each $z^* \in C^-\bs\cb{0}$. A closed convex function $f$ is proper if, and only if, there is $z^* \in C^-\bs\cb{0}$ such that $\vp_{f, z^*}$ is proper (in the usual sense of classical convex analysis), and this is the case if, and only if, the function $x \mapsto f\of{x} \oplus H\of{z^*}$ is proper. A standard separation argument shows
\[
\forall x \in X \colon f\of{x} = \bigcap_{z^* \in C^-\bs\cb{0}}\cb{z \in Z \mid \vp_{f, z^*}\of{x} + z^*\of{z} \leq 0}.
\]
With some effort, one can show that for a closed convex proper function $f \colon X \to \G^\triup$ it suffices to run the intersection in the above formula over the set of $z^* \in C^-\bs\cb{0}$ which generate a closed proper (and convex) scalarization $\vp_{f, z^*}$, see \cite{Schrage09Diss} and \cite[Corollary 3.34]{Schrage10-1ArX}.

\section{Directional derivatives and subdifferentials of $\G^\triup$-valued functions}
\label{SecDirDer}

\begin{definition}
\label{DefDirDer} Let $f \colon X \to \G^\triup$ be a convex function. The directional derivative of $f$ with respect to $z^*\in C^-\bs\cb{0}$ at $x_0 \in X$ in direction $x \in X$ is given by

\begin{equation}
\label{EqDirDerDef} 
f_{z^*}'\of{x_0, x} = \liminf_{t \downarrow 0} \frac{1}{t}\sqb{f\of{x_0 + tx} \istardif f\of{x_0}}
	= \bigcap_{s>0} \cl\bigcup_{0 < t < s} \frac{1}{t}\sqb{f\of{x_0 + tx} \istardif f\of{x_0}}.
\end{equation}
\end{definition}

If $f\of{x_0} = \emptyset$ then $f_{z^*}'\of{x_0, x} = Z$ for all $x \in X$. Therefore,
we can restrict the analysis to the case $x_0 \in \dom f$. The main tool will be the
directional difference quotient of $f$ at $x_0 \in X$ which is
defined to be the function $g_{z^*} \colon \R \times X \to \mathcal G^\triup$ given by
\[
g_{z^*}\of{t, x}  = \frac{1}{t}\sqb{f\of{x_0 + tx} \istardif f\of{x_0}}.
\]
The next lemma demonstrates the monotonicity of the difference quotient.

\begin{lemma}
\label{LemIncreasDiffQuot} Let $f \colon X \to \G^\triup$ be convex and $x_0, x \in X$. If $0 < t \leq s$ then
\begin{equation}
\label{EqIncreasDiffQuot1}
 H\of{z^*} \istardif g_{z^*}\of{s, -x} \supseteq H\of{z^*} \istardif g_{z^*}\of{t, -x}
    \quad \text{and} \quad g_{z^*}\of{t, x} \supseteq g_{z^*}\of{s, x}.
\end{equation}
If, additonally, $f\of{x_0} \oplus H\of{z^*} \not\in \cb{Z, \emptyset}$, then
\begin{equation}
\label{EqIncreasDiffQuot2}
H\of{z^*} \istardif g_{z^*}\of{t, -x} \supseteq g_{z^*}\of{t, x}.
\end{equation}
\end{lemma}

{\sc Proof.} Since $0 < \frac{t}{s} \leq 1$, $x_0 + tx = \frac{t}{s}\of{x_0 + s x} +
\frac{s-t}{s} x_0$ is well-defined. The convexity of $f$ produces
\[
f\of{x_0 + tx} \supseteq
    \frac{t}{s}f\of{x_0 + sx} + \frac{s-t}{s}f\of{x}.
\]
The rules (a), (c) and (b) of proposition \ref{PropCalcIstardif} produce
\begin{align*}
g_{z^*}\of{t, x}
    & \supseteq
    \frac{1}{t} \sqb{\of{\frac{t}{s}f\of{x_0 + sx} +
    \frac{s-t}{s}f\of{x_0}} \istardif f\of{x_0}} \\
    & \supseteq
    \frac{1}{t} \sqb{\frac{t}{s}\of{f\of{x_0 + sx} \istardif f\of{x_0}} +
    \frac{s-t}{s}\of{f\of{x_0} \istardif f\of{x_0}}} \\
    & \supseteq \frac{1}{s}\sqb{f\of{x_0 + sx} \istardif f\of{x_0}}
    \oplus H\of{z^*} = g_{z^*}\of{s, x}.
\end{align*}
Hence $g_{z^*}\of{t, x} \supseteq g_{z^*}\of{s, x}$. Replacing $x$ by $-x$ we obtain $g_{z^*}\of{t, -x} \supseteq g_{z^*}\of{s, -x}$. Proposition \ref{PropCalcIstardif}, (a) produces
\[
H\of{z^*} \istardif g_{z^*}\of{s, -x} \supseteq H\of{z^*} \istardif g_{z^*}\of{t, -x}.
\]

It remains to demonstrate the inequality \eqref{EqIncreasDiffQuot2}. Since $x_0
= \frac{1}{2}\of{x_0+tx} + \frac{1}{2}\of{x_0-tx}$ convexity of $f$ implies
\[
f\of{x_0} \supseteq \frac{1}{2}f\of{x_0+tx} + \frac{1}{2}f\of{x_0-tx},
\]
and from proposition \ref{PropCalcIstardif}, (a), (c) and $f\of{x_0} \istardif f\of{x_0}
= H\of{z^*}$ we obtain
\begin{align*}
H\of{z^*} & \supseteq
    \sqb{\frac{1}{2}f\of{x_0+tx} \oplus \frac{1}{2}f\of{x_0-tx}} \istardif f\of{x_0} \\
    & \supseteq \frac{1}{2}\sqb{f\of{x_0 + tx} \istardif f\of{x_0}} \oplus
        \frac{1}{2}\sqb{f\of{x_0 - tx} \istardif f\of{x_0}} \\
    & =t \cdot  \sqb{\frac{1}{2}g_{z^*}\of{t, x} \oplus \frac{1}{2}g_{z^*}\of{t, -x}}.
\end{align*}
Since $H\of{z^*} $ is a cone, the above relation can be divided by $t > 0$. Proposition \ref{PropCalcIstardif}, (g) yields
\[
H\of{z^*} \istardif g_{z^*}\of{t, -x} \supseteq g_{z^*}\of{t, x}.
\]

This completes the proof of the lemma. \pend

\begin{lemma}
\label{LemDirDerSublinear} Let $f \colon X \to \G^\triup$ be convex, $x_0 \in X$ and
$z^* \in C^-\bs\cb{0}$. Then

\begin{equation}\label{EqInfDirDer}
\forall x \in X \colon f'_{z^*}\of{x_0, x} =  \inf_{t>0}\frac{1}{t}\sqb{f\of{x_0 + tx} \istardif f\of{x_0}},
\end{equation}
and the function
\[
x \mapsto f'_{z^*}\of{x_0, x}
\]
is sublinear as a function from $X$ into $\mathcal G^\triup$. If $x_0 \in \dom f$, then $\dom f'_{z^*}\of{x_0,
\cdot} = \cone\of{\dom f - x_0}$. Moreover,
\[
f'_{z^*}\of{x_0, 0} =
    \left\{
        \begin{array}{ccc}
        H\of{z^*} & : & f\of{x_0} \oplus H\of{z^*} \not\in \cb{Z, \emptyset} \\
        Z & : & f\of{x_0} \oplus H\of{z^*} \in \cb{Z, \emptyset}
        \end{array}
    \right. .
\]
\end{lemma}

{\sc Proof.} The monotonicity of the difference quotient as proven in lemma \ref{LemIncreasDiffQuot} yields
\[
\bigcup_{0 < t < s} \frac{1}{t}\sqb{f\of{x_0 + tx} \istardif f\of{x_0}} = \bigcup_{t > 0} \frac{1}{t}\sqb{f\of{x_0 + tx} \istardif f\of{x_0}},
\]
and this implies \eqref{EqInfDirDer}. In turn,  \eqref{EqInfDirDer} immediately yields the positive homogeneity of $f'_{z^*}\of{x_0, \cdot}$. Take $x_1, x_2 \in X$. By convexity of $f$
\begin{align*}
f\of{x_0 + t\of{x_1 + x_2}} & =
    f\of{\frac{1}{2}\of{x_0 + 2tx_1} + \frac{1}{2}\of{x_0 + 2tx_2}} \\
    & \supseteq
    \frac{1}{2}f\of{x_0 + 2tx_1} + \frac{1}{2}f\of{x_0 + 2tx_2}.
\end{align*}
Proposition \ref{PropCalcIstardif}, (c) gives
\begin{multline}
\label{LemDirDerSublinear1}
\frac{1}{t}\sqb{f\of{x_0 + t\of{x_1 + x_2}} \istardif f\of{x_0}} \supseteq \\
    \frac{1}{2t}\sqb{f\of{x_0 + 2tx_1} \istardif f\of{x_0}}
    + \frac{1}{2t}\sqb{f\of{x_0 + 2tx_2} \istardif f\of{x_0}}.
\end{multline}
Fix $s>0$ and choose $t>0$ such that $2t \leq s$. Then, by the monotonicity of the
difference quotient (see \eqref{EqIncreasDiffQuot1} in Lemma \ref{LemIncreasDiffQuot})
\begin{multline*}
\frac{1}{2t}\sqb{f\of{x_0 + 2tx_1} \istardif f\of{x_0}}
    + \frac{1}{2t}\sqb{f\of{x_0 + 2tx_2} \istardif f\of{x_0}}  \supseteq \\
    \frac{1}{2t}\sqb{f\of{x_0 + 2tx_1} \istardif f\of{x_0}}
    + \frac{1}{s}\sqb{f\of{x_0 + sx_2} \istardif f\of{x_0}}.
\end{multline*}
This implies
\begin{multline*}
\cl \bigcup_{t > 0} \sqb{\frac{1}{2t}\sqb{f\of{x_0 + 2tx_1} \istardif f\of{x_0}}
    + \frac{1}{2t}\sqb{f\of{x_0 + 2tx_2} \istardif f\of{x_0}}}  \supseteq \\
\cl \bigcup_{0 < t \leq s} \frac{1}{2t}\sqb{f\of{x_0 + 2tx_1} \istardif f\of{x_0}}
    + \frac{1}{s}\sqb{f\of{x_0 + sx_2} \istardif f\of{x_0}}
\end{multline*}
for all $s > 0$ which produces, together with \eqref{LemDirDerSublinear1}, the desired
subadditivity.

We have $f'_{z^*}\of{x_0, x} = \emptyset$ if and only if
\[
\forall t > 0 \colon f\of{x_0 + tx} \istardif f\of{x_0} = \emptyset.
\]
Since, by assumption, $f\of{x_0} \neq \emptyset$, this is true if and only if $f\of{x_0 +
tx} = \emptyset$ for all $t>0$ (see definition of $\istardif$). This proves $\dom
f'_{z^*}\of{x_0, \cdot} = \cone\of{\dom f - x_0}$.

Finally, one easily checks the formula for $f'_{z^*}\of{x_0, 0}$. \pend

\begin{example}
A function $\vp \colon X \to \OLR$ can be identified with the function $f \colon X \to \G\of{\R, \R_+}$ defined by $f\of{x} = \vp\of{x} + \R_+$ where it is understood that $f\of{x} = \R$ whenever $\vp\of{x} = -\infty$ and $f\of{x} = \emptyset$ whenever $\vp\of{x} = +\infty$ (compare example \ref{ExExtReals}). Vice versa, if  $f \colon X \to \G\of{\R, \R_+}$ is given, one obtains a function $\vp \colon X \to \OLR$ by $\vp\of{x} = \inf\cb{r \in \R \mid r \in f\of{x}}$ with $\vp\of{x} = -\infty$ for $f\of{x} = \R$ and $\vp\of{x} = +\infty$ for $f\of{x} = \emptyset$. Obviously, $\vp$ is convex if, and only if, $f$ is convex. In this case, we obtain a directional derivative $f'$ (one and the same for all $z^* = s \in \of{\R_+}^-\bs\cb{0} = -\R_+\bs\cb{0}$)  for $f$ by means of definition \ref{DefDirDer}. We set for $x \in X$
\begin{equation}
\label{EqScalarDirDer}
\vp'\of{x_0, x} = \inf\cb{r \in \R \mid r \in f'\of{x_0, x}} ,
\end{equation}
and this definition is an extension of the classical definition for the directional derivative of proper extended real-valued functions. One can show (see \cite{HamelSchrage12}) that for convex $f$
\[
\forall x \in X \colon \vp'\of{x_0, x} = \inf\limits_{t>0}\frac{1}{t}\sqb{\vp\of{x_0+tx} \idif \vp\of{x_0}}.
\]
\end{example}

\begin{remark} \label{RemScalarDerivative} The function $\psi \colon X \to \OLR$ defined by
\[
\psi\of{x} = \inf\cb{-z^*\of{z} \mid z \in f'_{z^*}\of{x_0, x}}
\]
is the scalarization of $x \mapsto f'_{z^*}\of{x_0, x}$. Since $f'_{z^*}\of{x_0, \cdot}$ maps into $\mathcal G\of{H\of{z^*}}$, it is an exercise to prove that for all $x \in X$
\[
\psi\of{x}  = \vp'_{f, z^*}\of{x_0, x}
\]
and
\[
f'_{z^*}(x_0,x)=\cb{z \in Z \mid \psi\of{x}  + z^*\of{z} \leq 0}.
\]
These formulas basically tell us that the two operations 'scalarization' and 'taking the directional derivative' commute. The reader may compare \cite{Schrage09Diss}.
\end{remark}

\begin{example}
\label{ExPolyFunction1}
Let $X = Z = \R^2$ and $C = \R^2_+$. Consider $f \colon X \to \G\of{Z,C}$ defined by
\[
f\of{x} = \left\{
\begin{array}{ccc}
	\cb{z\in \R^2 \mid z_1\geq -x_1+x_2,\, z_2\geq -x_1-x_2,\, z_1+z_2 \geq x_1} & :  &x_1\geq 0 \\
	\emptyset & : & \text{otherwise}
\end{array}
\right.
\]
The function $f$ is convex since $\gr f$ is convex. Since $f\of{x}$ is given through three linear inequalities it can have at most two vertexes which depend linearly on $x$:
\[
V^1\of{x} = \of{\begin{array}{c}
  -x_1 + x_2 \\ 
  2x_1 - x_2\\ 
\end{array}} \quad \text{and} \quad
V^2\of{x} = \of{\begin{array}{c}
   2x_1 + x_2 \\ 
   -x_1 - x_2\\ 
\end{array}}.
\]
Moreover, $V^2\of{x} - V^1\of{x} = 3x_1\of{1, -1}^T$, and the recession cone of $f\of{x}$ is $\R^2_+$. This shows that for $w \in C^- = -\R^2_+$
\[
f\of{x} \oplus H\of{w} = 
\left\{
\begin{array}{ccc}
   V^1\of{x} + H\of{w} & : & w_1 \geq w_2 \\ 
  V^2\of{x} + H\of{w} & : & w_1 \leq w_2 \\ 
\end{array}
\right.
\]
for all $x \in \dom f$. Hence
\[
\frac{1}{t}\sqb{f\of{x + ty} \istardif f\of{x}} =
	\left\{
\begin{array}{ccc}
   V^1\of{y} + H\of{w} & : & w_1 \geq w_2 \\ 
  V^2\of{y} + H\of{w} & : & w_1 < w_2 \\ 
\end{array}
\right.
\]
whenever $x_1 > 0$ or $y_1 \geq 0$. Since this 'difference quotient' does not depend on $t$, it equals the directional derivative. We obtain for $x_1 > 0$ or $y_1 \geq 0$
\[
f'_w\of{x, y} =  \left\{
\begin{array}{ccc}		
		V^1\of{y} + H\of{w} & : & w_1\geq w_2 \\
		 V^2\of{y} + H\of{w}  & : & w_1 < w_2,
\end{array}
\right.
\]
and $f'_w\of{x, y} = \emptyset$ whenever $x_1 = 0$ and $y_1 < 0$. Since $z \in V^i\of{y} + H\of{w}$, if and only if, $w^T\of{z - V^i\of{y}} \leq 0$, $i = 1,2$,
\begin{align*}
V^1\of{y} + H\of{w} & = \cb{z \in \R^2 \mid y_1\of{w_1 - 2w_2} + y_2\of{-w_1 + w_2} + w^Tz \leq 0}, \\
V^2\of{y} + H\of{w} & = \cb{z \in \R^2 \mid y_1\of{-2w_1 + w_2} + y_2\of{-w_1 + w_2} + w^Tz \leq 0}.
\end{align*}
Finally, $f'_w\of{x, y} = Z$ for $x \not\in \dom f$, i.e. $x_1 < 0$.
\end{example}

The following result tells us when the directional derivative has only "finite" values. As usual, we denote by $\core M$ the algebraic interior of a set $M \subseteq X$.

\begin{theorem}
\label{ThmDirDerFinite} Let $f \colon X \to \G^\triup$ be convex and $x_0 \in
\core\of{\dom f}$. If $f$ is proper ($C$-proper), then there exists $z^* \in
C^-\bs\cb{0}$ ($z^* \in C^-\bs-C^-$) such that $f'_{z^*}\of{x_0, x} \not\in \cb{Z,
\emptyset}$ for all $x \in X$.
\end{theorem}

{\sc Proof.} Since $f$ is proper ($C$-proper), we have $f\of{x_0} \neq Z$. Hence there is
$z^* \in C^-\bs\cb{0}$ ($z^* \in C^-\bs-C^-$) such that $f\of{x_0} \oplus H\of{z^*} \neq
Z$ (a separation argument). This implies $f'_{z^*}\of{x_0, 0} = H\of{z^*}$.

We have to show that $f'_{z^*}\of{x_0, x} \neq Z$ for all $x \in X$. Since $x_0 \in
\core\of{\dom f}$ there is $t_0 > 0$ such that
\[
\forall t \in \sqb{0, t_0} \colon x_0 \pm tx \in \dom f,
\]
hence $\pm x \in \cone\of{\dom f - x_0} = \dom f'_{z^*}\of{x_0, \cdot}$ (see lemma \ref{LemDirDerSublinear}). Sublinearity of $f'_{z^*}\of{x_0, \cdot}$ implies
\[
H\of{z^*} = f'_{z^*}\of{x_0, 0} \supseteq f'_{z^*}\of{x_0, x} \oplus f'_{z^*}\of{x_0, -x} \neq \emptyset.
\]
This implies $f'_{z^*}\of{x_0, x}, f'_{z^*}\of{x_0, -x} \neq \emptyset$. In turn, this and the latter inclusion imply $f'_{z^*}\of{x_0, x}, f'_{z^*}\of{x_0, -x} \neq Z$.
\pend

\medskip Using "linear" minorants of the sublinear directional derivative we define
elements of the subdifferential.

\begin{definition}
\label{DefSubdiff} Let $f \colon X \to \G^\triup$ be convex, $x_0 \in X$ and $z^* \in C^-\bs\cb{0}$. The set
\[
\partial_{z^*}f\of{x_0} = \cb{x^* \in X^* \mid \forall x \in X \colon
    S_{\of{x^*, z^*}}\of{x} \supseteq f_{z^*}'\of{x_0, x}}
\]
is called the $z^*$-subdifferential of $f$ at $x_0$.
\end{definition}

Note that
\[
\forall s > 0 \colon s\partial_{z^*} f\of{x} = \partial_{sz^*}f\of{x}
\]
by virtue of  $f_{sz^*}'\of{x_0, x} =  f_{z^*}'\of{x_0, x}$. This relationship will be used in the next section to establish a link with adjoints of processes. The directional derivative $f_{z^*}'\of{x_0, \cdot}$ is improper, if, and only if, $\partial_{z^*} f\of{x_0} = \emptyset$. Elements of the $z^*$-subdifferential can also be characterized by the subdifferential inequality.

\begin{proposition}
\label{PropSubdiffInequality} Let $f \colon X \to \G^\triup$ be convex and $x_0 \in X$. The following statements are equivalent for $x^* \in X^*$, $z^* \in C^-\bs\cb{0}$:

(a) $\forall x \in X$: $S_{\of{x^*, z^*}}\of{x} \supseteq f'_{z^*}\of{x_0, x}$,

(b) $\forall x \in X$: $S_{\of{x^*, z^*}}\of{x - x_0} \supseteq f\of{x} \istardif
f\of{x_0}$.

If, additionally, $\dom f \neq \emptyset$ and $f\of{x_0} \oplus H\of{z^*} \neq Z$, then (a) and (b) are equivalent to

(c) $\forall x \in X$: $f\of{x_0} \oplus S_{\of{x^*, z^*}}\of{-x_0} \supseteq f\of{x}
\oplus S_{\of{x^*, z^*}}\of{-x}$ and

(d) $\forall x \in X$: $S_{\of{x^*, z^*}}\of{x} \istardif f\of{x} \supseteq S_{\of{x^*,
z^*}}\of{x_0} \istardif f\of{x_0}$.
\end{proposition}

{\sc Proof.}  (a) $\Rightarrow$ (b): Choosing $t=1$ in the definition of the directional
derivative and replacing $x$ by $x-x_0$ we obtain (b) from (a).

(b) $\Rightarrow$ (a): Use $x = x_0 + ty$ for $t>0$, $y \in X$ and $S_{\of{x^*,
z^*}}\of{ty} = tS_{\of{x^*, z^*}}\of{y}$.

Next, note that if $x_0 \not\in \dom f$ and $\dom f \neq \emptyset$, then there is no pair $\of{x^*, z^*}$ which satisfies either of the conditions (a), (b), (c), (d). So, assume $x_0 \in \dom f$ in the remaining part of the proof.

(b) $\Rightarrow$ (c): Since $x_0 \in \dom f$, the additional assumptions gives $f\of{x_0}
\oplus H\of{z^*} \not\in \cb{Z, \emptyset}$. Thus, we can apply the equality case in
proposition \ref{PropCalcIstardif}, (f) with $A = f\of{x}$, $B = f\of{x_0}$.

(c) $\Rightarrow$ (b): Because $x_0 \in \dom f$ and the additional assumption is satisfied, we can apply the equality case in
proposition \ref{PropCalcIstardif}, (d) with $A = S_{\of{x^*, z^*}}\of{x - x_0}$ and $B =
f\of{x_0}$.

(c) $\Leftrightarrow$ (d): This follows with the help of proposition
\ref{MinusSetInequality} and
\begin{multline*}
H\of{z^*} \istardif \of{f\of{x} \oplus S_{\of{x^*, z^*}}\of{-x}} =
    \cb{z \in Z \;\vert\; f\of{x} + S_{\of{x^*, z^*}}\of{-x} + z \subseteq H\of{z^*}} \\
    = S_{\of{x^*, z^*}}\of{x} \istardif f\of{x}
\end{multline*}
since $S_{\of{x^*, z^*}}\of{-x} = -x^*\of{x}\hat{z} + H\of{z^*}$ for
$z^*\of{\hat{z}}=-1$. \pend

\medskip To complete the picture, we establish the relationship with the subdifferential of the scalarizations in the next proposition. Note that because we do not a priori exclude the improper case we cannot just use the known scalar results.

Using \eqref{EqScalarDirDer} we define
\[
\partial \vp_{f, z^*}\of{x_0} = \cb{x^* \in X^* \mid \forall x \in X \colon x^*\of{x} \leq \vp_{f, z^*}'\of{x_0, x}}
\]
which, according to \cite[Proposition 5.5]{HamelSchrage12}, coincides with the definition given in \cite[Definition 5.4]{HamelSchrage12}. Now, if $f \colon X \to \G^\triup$ is convex, $x_0 \in X$ and $x^* \in X^*$, $z^* \in C^-\bs\cb{0}$, then
\begin{equation}
\label{EqSubdiffScalar}
\partial_{z^*}f\of{x_0}  = \partial \vp_{f, z^*}\of{x_0}.
\end{equation}
Indeed, this follows from remark \ref{RemScalarDerivative} and example \ref{ExConlinFunction}.



\medskip The first main result of the paper is the following set-valued
extension of the so-called max-formula.

\begin{theorem}
\label{ThmMaxFormula} Let $f \colon X \to  \G^\triup$ be a convex function, $x_0 \in \dom
f$ and $z^* \in C^-\bs\cb{0}$ such that the function $x \mapsto  f\of{x} \oplus H\of{z^*}$
is proper and the function $\vp_{f, z^*} \colon X \to \R\cup\cb{+\infty}$ is upper semi-continuous
at $x_0$. Then $\partial_{z^*}f\of{x_0} \neq \emptyset$ and it holds
\begin{equation}
\label{EqMaxFormula} \forall x \in X \colon f'_{z^*}\of{x_0, x} =
    \bigcap_{x^* \in \partial_{z^*}f\of{x_0}} S_{\of{x^*, z^*}}\of{x}.
\end{equation}
Moreover, for each $x \in X$ there exists $x^*_0 \in \partial_{z^*}f\of{x_0}$ such that
\begin{equation}
\label{EqMaxFormulaMax} f'_{z^*}\of{x_0, x} = S_{\of{x^*_0, z^*}}\of{x}.
\end{equation}
\end{theorem}

{\sc Proof.} We have $\dom \sqb{f\of{\cdot}\oplus H\of{z^*}} = \dom \vp_{f, z^*}$. Since the extended real-valued convex proper function $\vp_{f, z^*}$ is continuous at $x_0$, $\partial\vp_{f, z^*}\of{x_0} \neq \emptyset$ (see Theorem 2.4.9 in \cite{Zalinescu02})and
\[
\forall x \in X \colon \vp'_{f, z^*}\of{x_0, x}
 = \sup_{x^* \in \partial\vp_{f,z^*}\of{x_0}} x^*\of{x},
\]
and for each $x \in X$ there is $x^*_0 \in \partial\vp_{f, z^*}\of{x_0}$ such that
$\vp'_{f, z^*}\of{x_0, x} =  x^*\of{x}$.

The scalarization formulas for the directional derivative (see remark \ref{RemScalarDerivative})  produce the desired result.
\pend

\begin{remark}
The regularity assumption in the previous theorem is concerned with the scalarization
$\vp_{f, z^*}$ of the set-valued function $f$. This might not seem appropriate. However,
the assumption used seems to be the weakest possible. For example, it is implied by the
following interior point condition: there is $z_0 \in Z$ such that $\of{x_0, z_0} \in
\Int\of{\gr f}$. For this and further details about continuity concepts for set-valued functions, compare \cite{HeydeSchrage11R}, for example proposition 3.25 and theorem 4.4.
\end{remark}

Finally, we shall link the subdifferential with Fenchel conjugates for set-valued
functions as introduced in \cite{Hamel09} and \cite{Schrage09Diss}. Recall that the function
$-f^* \colon X^* \times C^-\bs\cb{0} \to \G^\triup$ defined by
\[
-f^*\of{x^*, z^*} =
    \inf_{x \in X} \sqb{f\of{x} + S_{\of{x^*, z^*}}\of{-x}} =
    \cl\bigcup_{x \in X}\sqb{f\of{x} + S_{\of{x^*, z^*}}\of{-x}}
\]
for $x^* \in X^*$ and $z^* \in C^-\bs\cb{0}$ has been called the (negative) Fenchel
conjugate of $f$ and the function $f^* \colon X^* \times C^-\bs\cb{0} \to \G^\triup$
defined by
\[
f^*\of{x^*, z^*} = \sup_{x \in X} \sqb{S_{\of{x^*, z^*}}\of{x} \istardif f\of{x}}
    = \bigcap_{x \in X} \sqb{S_{\of{x^*, z^*}}\of{x} \istardif f\of{x}}
\]
for $x^* \in X^*$ and $z^* \in C^-\bs\cb{0}$ has been called the (positive) Fenchel
conjugate of $f$.

\begin{corollary}
\label{CorSubdiffConjugate} For $x_0 \in \dom f$, $x^* \in X^*$ and $z^* \in
C^-\bs\cb{0}$ with $f\of{x_0} \oplus H\of{z^*} \neq Z$ the following statements are
equivalent:

(a) $x^* \in \partial_{z^*}f\of{x_0}$,

(b) $-f^*\of{x^*, z^*} = f\of{x_0} \oplus S_{\of{x^*, z^*}}\of{-x_0}$,

(c) $f^*\of{x^*, z^*} = S_{\of{x^*, z^*}}\of{x_0} \istardif f\of{x_0}$.
\end{corollary}

{\sc Proof.} Apply proposition \ref{PropSubdiffInequality} (a), (c) and (d). \pend

\section{$\mathcal G^\triup$-valued functions and convex processes}
\label{SecConProc}

A convex process $F$ defined on $X$ and with values in $Z$ is a set-valued map whose graph is a convex cone, and a closed convex process is a convex process whose graph is closed, compare, for example \cite[Definition 2.1.1]{AubinFrankowska90Book}. That is, convex processes and sublinear functions mapping into $\mathcal G\of{Z, \cb{0}}$ represent the same concept. Note that the choice $C = \cb{0}$ admits to include the case where no cone is available a priori.

An important concept in the theory of convex processes is the notion of the adjoint process. If $F \colon X \to \mathcal G\of{Z, \cb{0}}$ is a convex process, its adjoint is $F^\diamond \colon Z^* \to \mathcal P\of{X^*}$ defined by
\[
F^\diamond\of{z^*} = \cb{x^* \in X^* \mid \forall \of{x, z} \in \gr F \colon x^*\of{x} \leq z^*\of{z}}
\]
(see  \cite[Definition 2.5.1]{AubinFrankowska90Book}). This is, 
\[
\of{z^*, x^*} \in \gr F^\diamond \quad \Leftrightarrow \quad \of{x^*, -z^*} \in \of{\gr F}^-.
\]
The definition of $F^\diamond$ readily implies
\[
x^* \in F^\diamond\of{-z^*} \quad \Leftrightarrow \quad 
	\forall x \in X \colon F\of{x} \subseteq S_{\of{x^*, z^*}}\of{x}.
\]
Thus, $x^* \in F^\diamond\of{-z^*}$ if, and only if, $S_{\of{x^*, z^*}}$ is a conlinear minorant of $F$. The collection of linear minorants of a scalar or even vector-valued sublinear function is sometimes called the support set of the sublinear function (see, for example, \cite[p. 119]{Rubinov77}). Since in our framework the functions $S_{\of{x^*, z^*}} \colon X \to  \mathcal G\of{Z, \cb{0}}$ replace continuous linear operators, $\gr F^\diamond$ can be identified with the support set of $F$. Moreover, if the sublinear function $F \colon X \to \mathcal G^\triup$ is closed and proper, then it is the pointwise supremum of its conlinear minorants (compare \cite[Proposition 14]{Hamel09}), and this produces
\[
F\of{x} = \bigcap_{z^* \in C^-\bs\cb{0}}\sqb{\bigcap_{x^*\in F^\diamond\of{-z^*} }S_{\of{x^*, z^*}}\of{x}}.
\]


The following proposition summarizes the situation for a convex set-valued function.

\begin{proposition}
\label{PropSubdiffConProc}
Let $f \colon X \to  \G^\triup$ be convex, $x_0 \in \dom f$ and $F\of{x} = f'_{z^*}\of{x_0, x}$ for $x \in X$. Then
\[
F^\diamond\of{u^*} = \left\{
	\begin{array}{ccl}
	\emptyset & : & u^* \not\in\R_+\cb{-z^*} \\
	\of{\cone\of{\dom f - x_0}}^- & : & u^* = 0 \\
	s\partial_{z^*}f\of{x_0} & : & u^* = -sz^*, \; s>0
	\end{array}
	\right.
\]
\end{proposition}

{\sc Proof.}  By definition of the adjoint process,
\begin{align*}
F^\diamond\of{u^*} & = \cb{x^* \in X^* \mid \forall \of{x, z} \in \gr F \colon x^*\of{x} \leq u^*\of{z}} \\
	& =  \cb{x^* \in X^* \mid \forall x \in X \colon f'_{z^*}\of{x_0, x} \subseteq S_{\of{x*, -u^*}}\of{x}}
\end{align*}
Since $f'_{z^*}\of{x_0, x}$ and $S_{\of{x*, -u^*}}\of{x}$ are both closed half spaces with normals $z^*$ and $-u^*$, respectively, $F^\diamond\of{u^*} = \emptyset$ if $u^* \not\in\R_+\cb{-z^*}$. If $u^* = 0$, then $x^* \in F^\diamond\of{u^*}$ if, and only if,
\[
\forall x \in \dom f'_{z^*}\of{x_0, \cdot} \colon x^*\of{x} \leq 0.
 \]
Lemma \ref{LemDirDerSublinear} yields the second case. Finally, if $u^* = -sz^*$ for some $s>0$, then $x^* \in F^\diamond\of{u^*}$ if, and only if,
\[
\forall x \in X \colon f'_{z^*}\of{x_0, x} \subseteq S_{\of{x*, sz^*}}\of{x} = S_{\of{\frac{1}{s}x*, z^*}}\of{x},
\]
therefore $\frac{1}{s}x^*\in\partial_{z*}f\of{x_0}$. Definition \ref{DefSubdiff} produces the result in this case. \pend

\medskip The above proposition shows that the directional derivative and the subdifferential of a convex $\G^\triup$-function $f$ are related via convex process duality. In particular,
\[
\partial_{z^*}f\of{x_0} = \of{f'_{z^*}\of{x_0, \cdot}}^\diamond\of{-z^*}.
\]

In \cite[Proposition 1.37]{Mordukhovich06-1Book}, it is proven that the coderivative of a convex-graph multifunction $f \colon X \to \P\of{Z}$ at some point $\of{x_0, z_0} \in \gr f$ coincides with the function $\hat D^*f\of{x_0, z_0} \colon Z^* \to \P\of{X^*}$ given by
\[
\hat D^*f\of{x_0,z_0}\of{-z^*} = \cb{x^* \in X^* \mid \of{x^*,z^*} \in \of{\cone\of{\gr f - (x_0,z_0)}}^-},
\]
and it is pointed out ''that the coderivatives under consideration can be viewed as proper set-valued generalizations of the adjoint linear operator to the classical derivative at the point in question'' (\cite[p. 46]{Mordukhovich06-1Book}). This construction is closely related to our subdifferentials.

If $f \colon X \to \G^\triup$ is a convex function and $\of{x_0, z_0} \in\gr f$, then the normal cone of $\gr f$ at $\of{x_0, z_0} \in \gr f$ 
\[
\mathcal{N}_{\gr f}\of{x_0,z_0} = \of{\cone\of{\gr f - \of{x_0, z_0}}}^-
\]
consists exactly of those elements $\of{x^*, z^* } \in X^* \times C^-$ with $z^*=0$ and $x^*\in\of{\cone\of{\dom f-x_0}}^-$, or $z^*\in C^-\setminus\cb{0}$,
$-z^*(z_0)=\vp_{f,z^*}(x_0)$ and $S_{\of{x^*,z^*}}(x_0)\istardif f(x_0)=f ^*(x^*,z^*)$. From corollary 3.11 we may conclude
\[
\hat D^*f\of{x_0, z_0}\of{-z^*} =
\left\{
	\begin{array}{ccc}
 	\partial_{z^*}f\of{x_0} & : & z^*\in C^-\bs\cb{0} \, \wedge \, -z^*\of{z_0} = \vp_{f,z^*}\of{x_0}\\
 	\of{\cone\of{\dom f-x_0}}^- & : & z^*=0 \\
 	\emptyset & : &  \text{otherwise} 
	\end{array}
\right.
\]
Hence for all $z^*\in C^-\setminus\cb{0}$
\[
\hat D^*f\of{x_0, z_0}\of{-sz^*} = 
\left\{
	\begin{array}{ccc}
	\of{f'_{z^*}\of{x_0,\cdot}}^\diamond\of{-sz^*} & : & s \geq 0 \; \text{and} \; -sz^*\of{z_0} = s\vp_{f,z^*}\of{x_0}\\
	\emptyset & : & \text{otherwise} 
	\end{array}
\right.
\]
where the convention $0\of{\pm\infty} =  0$ is essential. 

Finally, it may happen that the Mordukhovich coderivative is empty whereas the $z^*$-subdifferential is non-empty. For example, if $f \colon \R \to \G\of{\R^2,\R_+}$ given by $f\of{x} = \cb{\of{z_1, z_2} \in \R^2 \mid z_2 \geq \frac{1}{z_1}}$ for all $x \in \R$ and $z^*=\of{-1,0}$, then $\partial_{z^*}f\of{x} = \cb{0}$, and for all $z \in f\of{x} $ we have $-z^*\of{z} > \vp_{z^*}\of{x} = 0$. 

The relationship between Mordukhovich's coderivative and our directional derivative/subdifferential via convex process duality may be roughly summarized as follows. At least for convex set-valued functions, the set-valued directional derivative provides the 'primal concept' whose adjoint is a slight extension of Mordukhovich's coderivative.\footnote{This provides some evidence for the statement 'Primal is primal and always possible.' repeatedly given by J.-P. Penot at the Spring School on Variational Analysis, Paseky na Jizerou, 2003.}

\section{$\G^\triup$-valued optimization problems}

We are interested in the problem
\[
\tag{P} \mbox{minimize} \quad f \quad \mbox{subject to} \quad  x \in X
\]
where $f$ is a $\G^\triup$-valued function. The minimization is understood as looking for
the infimum in $\G^\triup$, that is $\inf_{x \in X} f\of{x} = \cl\co\cup_{x \in
X}f\of{x}$, and subsets of $X$ in which this infimum is attained. This approach, initiated in \cite{HeydeLoehne11}, is different from most other approaches in set optimization, see for example \cite[Definition 14.2]{Jahn04}, \cite{HernandezRodriguezMarin07-1}, \cite{HernandezRodriguezMarin07-2} and the references therein. Note that the convex hull in the definition of the infimum above can be dropped if $f$ is convex.

\begin{definition}
\label{DefInfimizer} Let $f \colon X \to \G^\triup$ be a function. A set $M \subseteq X$
is called an infimizer of $f$ if
\[
\inf_{x \in M} f\of{x} = \inf_{x \in X} f\of{x}.
\]
\end{definition}

Note that $\dom f$ is always an infimizer of $f$ which amounts to the necessity to
introduce further requirements to a solution of (P). In particular, the values $f\of{x}$
for $x \in M$ should satisfy additional conditions, e.g. be minimal in some sense. See
\cite{HeydeLoehne11}, \cite{Loehne11Book} for further motivations and corresponding concepts.

\begin{example}
\label{ExPolyFunction2} Let us again consider the function $f$ given in example \ref{ExPolyFunction1} and define the two sets
\[
M = \cb{x \in \R^2 \mid x_1 = 0} \quad \text{and} \quad N = \dom f\bs M.
\]
Then, $\inf f = \cb{z \in \R^2 \mid z_1 + z_2 \geq 0}$ and both $M$ and $N$ are infimizers. Indeed,
\begin{multline*}
f\sqb{M} = \bigcup_{x_2 \in \R}\cb{z \in \R^2 \mid z_1 \geq x_2, \; z_2 \geq -x_2, \; z_1 + z_2 \geq 0}   
	= \cb{z \in \R^2 \mid z_1 + z_2 \geq 0},
\end{multline*}
and by definition of $f$, $f\of{x} \subseteq \cb{z \in \R^2 \mid z_1 + z_2 \geq 0}$ for all $x \in \R^2$. The infimizer property for $N$ is a consequence of the closure operation which is part of the infimum in $\G\of{C}$.
\end{example}

In this note, we introduce another solution concept for problem (P) and provide
characterizations in terms of the directional derivative and the subdifferential
introduced above.

\begin{definition}
Let $f \colon X \to \G^\triup$ be a function and $z^* \in C^-\bs\cb{0}$. A point $x_0 \in X$ is called a
$z^*$-solution of (P), if 
\[
f\of{x_0} \oplus H\of{z^*} = \sqb{\inf_{x \in X} f\of{x}} \oplus H\of{z^*}.
\]
A set $M \subseteq X$ is called a $C^-$-solution of (P) if $M$ is an infimizer of (P) and each $m \in
M$ is a $z^*$-solution for some $z^* \in C^-\bs\cb{0}$.
\end{definition}

 The concept of $z^*$-solutions clearly generalizes the solution of a real-valued
minimization problem as well as weak solutions in vector optimization (see \cite{Jahn04}, theorem 5.13). Moreover, it is also clear that $x_0$ is a $z^*$-solution of (P) if, and only if, it is a solution of the scalarized problem
\[
\mbox{minimize} \quad \vp_{f, z^*} \quad \mbox{subject to} \quad  x \in X.
\]

\begin{example}
\label{ExPolyFunction3} We continue the discussion of example \ref{ExPolyFunction1}. One may see that each $f\of{x}$ for $x \in \dom f$ is minimal with respect to $\supseteq$, i.e. $f\of{x} \supseteq f\of{y}$ implies $f\of{x} = f\of{y}$ (and even $x = y$). However, the only $w$-solutions for some $w \in C^-\bs\cb{0} = -\R^2_+\bs\cb{0}$ with a 'finite' value of $f$ are the elements of $M$, i.e. $x \in \R^2$ with $x_1 = 0$ (see example \ref{ExPolyFunction2}), and one has to choose $w$ such that $w_1 = w_2 < 0$. In particular, this shows that there are infimizers which do not include a single $w$-solution, and, moreover, that there are infimizers which are solutions in the sense of \cite{HeydeLoehne11}, \cite[Definition 2.8]{Loehne11Book}, but do not include a single $w$-solution. In fact, the set $N$ (example \ref{ExPolyFunction2}) is such an infimizer.
\end{example}

\begin{proposition}
\label{PropZStarSolutions} Let $f \colon X \to \G^\triup$ be convex and $z^*\in
C^-\bs\cb{0}$. The following statements are equivalent for $x_0 \in \dom f$:
\\
(a) $x_0$ is a $z^*$-solution of (P);
\\
(b) $\dom f = \emptyset$, or $f\of{x_0} \oplus H\of{z^*} = Z$, or  $H\of{z^*} \supseteq f'_{z^*}\of{x_0, x}$ for all $x \in
X$;
\\
(c) $\dom f = \emptyset$, or $f\of{x_0} \oplus H\of{z^*} = Z$, or $0 \in \partial_{z^*}f\of{x_0}$.
\end{proposition}

{\sc Proof.} Since the improper cases are trivial, let us assume $f\of{x_0} \oplus H\of{z^*} \not\in \cb{Z, \emptyset}$. 

(a) $\Rightarrow$ (b): Since $x_0$ is a $z^*$-solution of (P) we have
$f\of{x_0}  \oplus H\of{z^*} = \sqb{\inf_{x \in X} f\of{x}} \oplus H\of{z^*} \neq Z$. Take $x \in X$,
$t>0$ and $z \in f\of{x_0 + tx} \istardif f\of{x_0}$. Then
\[
f\of{x_0 + tx} \oplus H\of{z^*} + z \subseteq f\of{x_0} \oplus H\of{z^*} + z \subseteq f\of{x_0 + tx} \oplus H\of{z^*}.
\]
This is only possible if either $ f\of{x_0 + tx} \oplus H\of{z^*} = Z$ and hence $f\of{x_0} \oplus H\of{z^*} = Z$, or $z \in H\of{z^*}$. The former case is now excluded, so we may conclude (b).

(b) $\Rightarrow$ (a): (b) implies (choose $t=1$ and use the definition of the infimum)
\[
\forall x \in X \colon H\of{z^*} \supseteq f\of{x_0 + tx} \istardif f\of{x_0},
\]
hence
\[
\forall y \in X \colon H\of{z^*} \supseteq f\of{y} \istardif f\of{x_0}.
\]
Proposition \ref{PropCalcIstardif}, (g) implies
\[
\forall y \in X \colon f\of{y} \subseteq f\of{x_0} \oplus H\of{z^*} .
\]

The equivalence of (b) and (c) is proposition \ref{PropSubdiffInequality} for $x^* = 0$.
\pend

\medskip In most cases, a single $z^*$-solution does not characterize $\inf\limits_{x\in X}f\of{x}$ too
well. The final task in this paper is to characterize infimizers in terms of directional
derivatives and subdifferentials. In order to do so, we need one more new concept. Its introduction is motivated by the fact that an infimizer does not necessarily consist of $z^*$-solutions. On the other hand, the singleton $M = \cb{x_0}$ with $f\of{x_0} \oplus H\of{z^*} \neq Z$ is an infimizer if, and only if, it is a $z^*$-solution. One can understand the following procedure as a way to reduce infimizers to singletons.

\begin{definition}\label{def:CanExt}
Let $f \colon X \to \G^\triup$ be a function and $M \subseteq X$. We define the
inf-translation of $f$ by $M$ to be the function $\hat f\of{\cdot; M} \colon X \to
\G^\triup$ given by
\begin{equation}\label{EqSetTranslation}
\hat f\of{x; M} = \inf\limits_{m \in M}f\of{m + x} = \cl\co\bigcup_{m \in M}f\of{m+x}.
\end{equation}
\end{definition}

The function $\hat f\of{\cdot; M}$ is nothing else than the canonical extension of $f$ at
$M + x$ as defined in \cite{HeydeLoehne11}. 

\begin{example}
\label{ExPolyFunction4} We compute $\hat f\of{x; M}$ for $f$ and $M$ as defined in example \ref{ExPolyFunction1} and \ref{ExPolyFunction2}. The definitions of $f$ and $M$ imply $\hat f\of{x; M} = \emptyset$ for $x \not\in \dom f$. If $x \in \dom f$, then
\begin{align*}
\hat f\of{x; M} & = \cl\bigcup_{m_1 = 0, \, m_2 \in \R} f\of{m + x} \\
	& = \bigcup_{m_2 \in \R} 
	\cb{z \in \R^2 \mid z_1 \geq -x_1 + x_2 + m_2, \; z_2 \geq -x_1 - x_2 - m_2, \; z_1 + z_2 \geq x_1} \\
	& = \cb{z \in \R^2 \mid z_1 + z_2 \geq x_1}.
\end{align*}
\end{example}

A few elementary properties of the inf-translation are
collected in the following lemma.

\begin{lemma}\label{LemInfShiftProps}
Let $M \subseteq X$ be non-empty and $f \colon X \to \G^\triup$ a function. Then
\\
(a) if $M \subseteq N \subseteq X$ then $\hat f\of{x; M} \subseteq \hat f\of{x; N}$ for
all $x \in X$;
\\
(b) $\inf_{x \in X} f\of{x} = \inf_{x \in X} \hat f\of{x; M}$;
\\
(c) if $f$  and $M$ are convex, so is $\hat f\of{\cdot; M} \colon X \to \G^\triup$, and in this case $\hat f\of{x; M} = \cl\bigcup_{m \in M}f\of{m+x}$.
\end{lemma}

{\sc Proof.} (a) Immediate from the definition of $\hat f\of{\cdot; M}$.

(b) By definition of $\hat f\of{\cdot; M}$,
\begin{align*}
\inf_{x \in X} \hat f\of{x; M} & = \inf_{x \in X}\inf_{m \in M} f\of{m + x}
     = \inf_{m \in M}\inf_{x \in X} f\of{m + x} \\
    & = \inf_{m \in M}\sqb{\inf_{m + x \in X} f\of{m + x}}
    = \inf_{m \in M} \sqb{\inf_{y \in X}f\of{y}} = \inf_{x \in X}f\of{x}.
\end{align*}

(c) Take $t \in \of{0,1}$, $x_1, x_2 \in X$. Since $M$ is convex, we have $M = tM +
\of{1-t}M$. This and the convexity of $f$ yield
\begin{align*}
\hat f\of{tx_1+\of{1-t}x_2; M} & = \inf_{m \in M}f\of{tx_1+\of{1-t}x_2 + m} \\
    & = \inf_{m_1, m_2 \in M}f\of{t\of{x_1 + m_1} + \of{1-t}\of{x_2+m_2}} \\
    & \supseteq \inf_{m_1, m_2 \in M}tf\of{x_1 + m_1} \oplus \of{1-t}f\of{x_2+m_2} \\
    & = t\hat f\of{x; M} \oplus \of{1-t}\hat f\of{x_2; M}.
\end{align*}

This completes the proof of the lemma. \pend

\begin{proposition}
\label{PropInfShiftChar} Let $f \colon X \to \G^\triup$ be a convex function and
$\emptyset \neq M \subseteq \dom f$. The following statements are equivalent:
\\
(a) $M$ is an infimizer for $f$;
\\
(b) $\cb{0} \subseteq X$ is an infimizer for $\hat f\of{\cdot; M}$;
\\
(c) $\cb{0}$ is an infimizer for $\hat f\of{\cdot; \co M}$ and $\hat f\of{0; M} = \hat
f\of{0; \co M}$.
\end{proposition}

{\sc Proof.} The equivalence of (a) and (b) is immediate from $\hat f\of{0; M} = \inf_{m
\in M}f\of{m}$ and lemma \ref{LemInfShiftProps}, (b). The equivalence of (a) and (c)
follows from $\hat f\of{0; \co M} = \inf_{m \in \co M}f\of{m}$ and lemma \ref{LemInfShiftProps}
(b). \pend

\medskip It should be apparent that we need to consider $\hat f\of{\cdot; \co M}$: Since
we want to characterize infimizers via directional derivatives and subdifferentials, a
convex function is needed, and $\hat f\of{\cdot; M}$ is not convex in general even if $f$
is convex. Obviously, an infimizer need not be a convex set; on the contrary, sometimes
one prefers a nonconvex one, for example a collection of vertexes of a polyhedral set
instead of a face.

\begin{theorem}
\label{ThmOptimality} Let $f \colon X \to \G^\triup$ be a proper convex function with
\[
I\of{f} = \inf_{x \in X}f\of{x} \neq Z.
\]
Let $\Gamma^-\of{f} = \cb{z^* \in C^-\bs\cb{0} \mid I\of{f}\oplus H\of{z^*} \neq Z}$. The
following statements are equivalent:
\\
(a) $M$ is an infimizer for $f$;
\\
(b) $\hat f\of{0; M} = \hat f\of{0; \co M}$ and
\[
\forall z^* \in \Gamma^-\of{f}, \; \forall x \in X \colon H\of{z^*} \supseteq \hat
f'_{z^*}\of{\cdot; \co M}\of{0, x};
\]
(c) $\hat f\of{0; M} = \hat f\of{0; \co M}$ and
\[
0 \in \bigcap_{z^* \in \Gamma^-\of{f}}\partial_{z^*}\hat f\of{\cdot; \co M}\of{0}.
\]
\end{theorem}

{\sc Proof.} Since $\cb{0}$ is a singleton infimizer of the function $x \mapsto \hat f\of{x; M}$, $0 \in X$ is a $z^*$-solution of (P) for  $f\of{\cdot; M}$ for each $z^* \in \Gamma^-\of{f}$. Now, the result follows from proposition \ref{PropZStarSolutions} and proposition
\ref{PropInfShiftChar} .  \pend

\begin{example}
\label{ExPolyFunction5} We know that $M$ from  example \ref{ExPolyFunction2} is an infimizer of $f$ from example \ref{ExPolyFunction1}. Moreover,
\[
\hat f\of{x; M} = \inf_{x \in \R^2} f\of{x} =  \cb{z \in \R^2 \mid z_1 + z_2 \geq x_1},
\] 
for $x \in \dom f$, see example \ref{ExPolyFunction4}. Since $\hat f\of{x; M}$ is a closed half space with normal $\of{-1, -1}^T$ for all $x \in \dom f$ we obtain
\[
\of{\hat f\of{\cdot; M}}'_w\of{0, y} = 
	\left\{
	\begin{array}{ccc}
	\emptyset & : & y_1 < 0 \\
	\cb{z \in Z \mid w^Tz \leq w_1y_1} & : & w_1 = w_2, \; y_1 \geq 0 \\
	Z & : & w_1 \neq w_2, \; y_1 \geq 0
	\end{array}
	\right.
\]
This relationship illustrate (b) in theorem \ref{ThmOptimality} since for $w \in C^-\bs\cb{0} = -\R^2_+\bs\cb{0}$ with $w_1 = w_2$ we have $w_1y_1 \leq 0$ and so $\cb{z \in Z \mid w^Tz \leq w_1y_1} \subseteq H\of{w}$ whenever $y_1 \geq 0$. Obviously, for $w \in -\R^2_+\bs\cb{0}$
\[
 \inf_{x \in X}f\of{x} \oplus H\of{w} \neq \R^2 \quad \Leftrightarrow \quad w_1 = w_2
\]
since $ \inf_{x \in X}f\of{x} = \cb{z \in Z \mid z_1 + z_2 \geq 0}$ as shown in example  \ref{ExPolyFunction2}.
\end{example}

{\bf Acknowledgement.} The function $f$ for the 'running example' discussed in example \ref{ExPolyFunction1}, \ref{ExPolyFunction2}, \ref{ExPolyFunction3}, \ref{ExPolyFunction4} and \ref{ExPolyFunction5} has been constructed by Frank Heyde.

\end{document}